\documentclass{llncs}
\usepackage{epsfig,graphicx,amsfonts,amssymb,algorithm2e}

\def\bsquareforqed{\rule{0.6em}{0.6em}}
\def\bqed{\ifmmode\bsquareforqed\else{\unskip\nobreak\hfil
\penalty50\hskip1em\null\nobreak\hfil\bsquareforqed
\parfillskip=0pt\finalhyphendemerits=0\endgraf}\fi}
\spnewtheorem{observation}{Observation}{\bfseries}{\itshape}
\spnewtheorem{clm}{Claim}{\itshape}{\rmfamily}
\newcommand{\ds}{\displaystyle}  
\newcommand{\abs}[1]{| #1 |}

\newcommand{\cub}{\mbox{\textnormal{cub}}}
\newcommand{\boxi}{\mbox{\textnormal{box}}}
\newcommand{\mod}{\hspace{1mm}\mbox{\textnormal{mod}\hspace{1mm}}}
\newcommand{\ceil}[1]{\left\lceil #1 \right\rceil}
\newcommand{\floor}[1]{\left\lfloor #1 \right\rfloor}

\newcommand{\claw}{\psi(G)}
\newcommand{\clat}{\ceil{\log_2\alpha}}
\newcommand{\clct}{\ceil{\log_2\claw}}
\newcommand{\clcta}{\ceil{\log_2\psi}}

\newcommand{\loc}{\mathcal{C}}
\newcommand{\reals}{\mathbb{R}}
\newcommand{\integers}{\mathbb{Z}}

\begin{document}
\title{Cubicity of Interval Graphs and the Claw Number}
\author{Abhijin Adiga\thanks{Indian Institute of Science, Dept. of Computer
Science and Automation, Bangalore 560012, India. 
email: abhijin@csa.iisc.ernet.in}\and L. Sunil
Chandran\thanks{\textbf{(Corresponding Author)} Indian Institute of Science, Dept. of Computer
Science and Automation, Bangalore 560012, India. email: sunil@csa.iisc.ernet.in}}
\institute{}
\date{}
\maketitle
\begin{abstract}
Let $G(V,E)$ be a simple, undirected graph where $V$ is the set of
vertices and $E$ is the set of edges. A $b$-dimensional cube is a
Cartesian product $I_1\times I_2\times\cdots\times I_b$, where each $I_i$
is a closed interval of unit length on the real line. The \emph{cubicity}
of $G$, denoted by $\cub(G)$ is the minimum positive integer $b$ such
that the vertices in $G$ can be mapped to axis parallel $b$-dimensional
cubes in such a way that two vertices are adjacent in $G$ if and only
if their assigned cubes intersect. An interval graph is a graph that
can be represented as the intersection of intervals on the real line -
i.e., the vertices of an interval graph can be mapped to intervals on
the real line such that two vertices are adjacent if and only if their
corresponding intervals overlap. Suppose $S(m)$ denotes a star graph
on $m+1$ nodes. We define \emph{claw number} $\claw$ of the graph to
be the largest positive integer $m$ such that $S(m)$ is an induced
subgraph of $G$. It can be easily shown that the cubicity of any graph
is at least $\clct$.

In this paper, we show that, for an interval graph $G$
$\clct\le\cub(G)\le\clct+2$. It is not clear whether the upper bound of
 $\clct+2$ is tight: Till now we are unable to find any interval graph
with $\cub(G)>\clct$. We also show that, for an interval graph $G$,
 $\cub(G)\le\clat$, where $\alpha$ is the independence number of $G$.
Therefore, in the special case of $\claw=\alpha$, $\cub(G)$ is exactly
$\clat$.

The concept of cubicity can be generalized by considering boxes instead of
cubes. A $b$-dimensional box is a
Cartesian product $I_1\times I_2\times\cdots\times I_b$, where each $I_i$
is a closed interval on the real line.
The \emph{boxicity} of a graph, denoted $\boxi(G)$, is the minimum $k$ such
that $G$ is the intersection graph of $k$-dimensional boxes. It is clear
that $\boxi(G)\le\cub(G)$. From the above result, it follows that
for any graph $G$, $\cub(G)\le\boxi(G)\clat$.

\end{abstract}

\section{Introduction}
Let $G(V,E)$ be a simple, undirected graph where $V$ is the set
of vertices and $E$ is the set of edges. A $b$-dimensional box
is a Cartesian product $R_1\times R_2\times\cdots\times R_b$,
where each $R_i$ is a closed interval on the real line. When each
interval has unit length, we will call such a box a $b$-dimensional
cube. The \emph{cubicity} (respectively boxicity) of $G$, denoted by $\cub(G)$
($\boxi(G)$), is the minimum positive integer $b$ such that the vertices
in $G$ can be mapped to axis parallel $b$-dimensional cubes (boxes)
in such a way that two vertices are adjacent in $G$ if and only if
their assigned cubes (boxes) intersect. Cubicity and boxicity were
introduced by Roberts in \cite{recentProgressesInCombRoberts}. Yannakakis
\cite{complexityPartialOrderDimnYannakakis} proved that it is NP-complete
to determine if the cubicity of a graph is at most 3. It was shown by
Cozzens \cite{phdThesisCozzens} that computing the boxicity of a graph is
NP-hard. This was strengthened by Kratochvil
\cite{specialPlanarSatisfiabilityProbNPKratochvil} who showed that deciding
whether boxicity of a graph is at most 2 itself is NP-complete.

Roberts \cite{recentProgressesInCombRoberts} showed
 that for any graph $G$, $\cub(G)\le \floor{2n/3}$ and
 $\boxi(G)\le\floor{n/2}$. The cube
representation of special classes of graphs like hypercubes,
 co-bipartite and complete multipartite graphs were investigated in
\cite{cubicityHypercubeSunilNaveen,complexityPartialOrderDimnYannakakis,recentProgressesInCombRoberts}.
 Scheinerman \cite{phdThesisScheinerman} showed that the
 boxicity of outer planar graphs is at most 2. Thomassen
\cite{intervalRepPlanarGraphsThomassen} proved that the boxicity of
planar graphs is at most 3. In \cite{computingBoxicityCozzensRoberts},
 Cozzens and Roberts studied the boxicity of split graphs.
It is interesting to note that coloring problems on low boxicity
graphs were considered as early as 1948 \cite{problem56Bielecki}.
Kostochka \cite{coloringIntersectionGraphsKostochka} provides an
 extensive survey on colouring problems of intersection graphs. In
 \cite{findingImaiAsano,optimalPackingFowlerPatersonTanimoto} the
complexity of finding the maximum independent set in bounded boxicity
graphs is considered. In \cite{phdThesisHavel,onTheSphCubGraphsFishburn}
cubicity has been studied in comparison with sphericity. Some other related
references are
\cite{gridIntBoxicityBellantoniEtal,boxicityMaxDegreeSunilNaveenFrancis,boxicityTreewidthSunilNaveen,intNumBoxDigraphs,rectNumHypercubesChangeWest,forbiddenTrotter,posetboxicitytrotterwest,sphericityCubicityEdgeCliqueQuint}.

In this paper, we consider the cubicity of \emph{interval graphs}.
Graphs with boxicity at most 1 are precisely the well-studied class of
interval graphs. A graph is an interval graph if and only if its vertices
can be mapped to intervals on the real line such that two vertices are
adjacent if and only if their corresponding intervals overlap. From the
definition of boxicity and cubicity, it is easy to see that any cube
representation of a graph will also serve as a box representation. Hence,
$\boxi(G)\le\cub(G)$. Therefore, it is indeed interesting to ask the
following question: what is the cubicity of a graph whose boxicity is 1?

Chandran and Mathew \cite{upperboundCubicityBoxicitySunilAshik}
showed that cubicity of an interval graph is at most
$\ceil{\log_2|V|}$. This was later improved to $\ceil{\log_2\Delta}+4$ in
\cite{cubicityIntervalGraphSunilFrancisNaveen}, where $\Delta$ is the
maximum degree of $G$. We improve this bound further. To state our
result, we first introduce a parameter called \emph{claw number} of a
graph. Recall that a star graph on $n$ vertices is the complete bipartite
graph $K_{1,n-1}$. We denote it by $S(n-1)$. 
\begin{definition}
The \textbf{claw number} $\claw$ of a graph $G$ is the largest positive
integer $m$ such that $S(m)$ is an induced subgraph of $G$. 
\end{definition}
Our result is as follows:
\begin{theorem}\label{thm:mainThm}
Let $G$ be an interval graph with claw number $\psi$.
\[
\clcta\le\cub(G)\le\clcta+2.
\]
\end{theorem}
It is not clear whether the upper bound of $\clcta+2$ is tight. We
have not been able to find any interval graph with cubicity greater
than $\clcta$. By slightly modifying the proof of Theorem
\ref{thm:mainThm}, we can also show that for any interval graph $G$,
$\cub(G)\le\clat$, where $\alpha$ is the independence number of $G$. Thus,
for the special case of $\psi=\alpha$, $\cub(G)$ is exactly
$\clat$. This in turn allows us to infer that for any graph $G$,
$\cub(G)\le\boxi(G)\clat$ (See the end of Section \ref{sec:proof}).


\subsection{Some Basic Properties and Results}
In this section, we mention some useful properties and results regarding
interval graphs and cubicity. 
A restricted form of interval graphs, that allow only intervals of unit
length are called \emph{indifference graphs}. They are also known as
\emph{unit interval graphs} or \emph{proper interval graphs}. We provide
an alternate definition which we make use of in later sections.
\begin{definition}{\textbf{Indifference graph:}} \label{defn:indiff}
 A graph $G(V,E)$ is an indifference graph if and only if there exists
 a function $\Pi:V\longrightarrow\reals$ such that for two distinct
 vertices $u$ and $v$, $u$ and $v$ are adjacent if and only if
$\abs{\Pi(u)~-~\Pi(v)}\le t$, for some fixed positive real number $t$.
\end{definition}
It is easy to see that a graph has cubicity 1 if and only if it is
an indifference graph.
\begin{property}{(See Golumbic
\textnormal{\cite{algGraphTheoryPerfectGraphsGolumbic}} for a proof.)}
\label{prop:consec1}
A graph $G$ is an interval graph if and only if its maximal cliques can
be linearly ordered such that for every vertex $u$ the maximal cliques
containing $u$ occur consecutively.
\end{property}
For a graph $G(V,E)$, let $G_i(V,E_i)$, $i\in\{1,2,\ldots,k\}$ be
such that $E=E_1\cap E_2\cap\cdots\cap E_k$. Then we say that $G$
is the \emph{intersection} of $G_i$'s $1\le i\le k$ and denote it as
$G=\ds\bigcap_{i=1}^{k}G_i$. Cubicity (Boxicity respectively) can be
stated in terms of intersection of indifference graphs (interval graphs)
as follows:
\begin{lemma}{\textnormal{Roberts
\cite{recentProgressesInCombRoberts}}}\label{lem:intcubbox}
The cubicity (boxicity) of a graph $G$ is the minimum positive integer $b$
such that $G$ is the intersection of $b$ indifference graphs (interval
graphs). Moreover, if $G=\bigcap_{i=0}^{m-1}G_i$, for some graphs $G_i$,
then, $\cub(G)\le\sum_{i=0}^{m-1}\cub(G_i)$ and
$\boxi(G)\le\sum_{i=0}^{m-1}\boxi(G_i)$.  
\end{lemma}
The following result is easy to prove.
\begin{lemma}\label{lem:inducedGraphCub}
Suppose $H$ is an induced subgraph of $G$, then $\cub(G)\ge\cub(H)$.
\end{lemma}

\section{Proof of Theorem \ref{thm:mainThm}}\label{sec:proof}
The lower bound is easy to see and is as follows. Since
the claw number of $G$ is $\psi$, it has  an induced
subgraph $S(\psi)$ and $\cub(S(\psi))=\clcta$ (See Roberts
\cite{recentProgressesInCombRoberts}). By Lemma \ref{lem:inducedGraphCub},
 $\cub(G)\ge\cub(S(\psi))=\clcta$. 

Our aim is to construct $\clcta+2$ indifference graphs and show that
$G$ is the intersection of these graphs, thereby proving the upper
bound. First, we describe a vertex numbering which is essential for
the construction of the indifference graphs.

\subsection{Vertex Labelling and the Primary Maximum Independent
Set}\label{sec:indSet}
Let $G(V,E)$ be an interval graph. Let $\loc:C_0, C_1,\ldots, C_{k-1}$
correspond to a linear ordering of maximal cliques satisfying Property
\ref{prop:consec1}, where $C_i$ corresponds to the set of vertices in the $i$th
maximal clique. For a vertex $u$, let $c_u=\{i|u\in C_i\}$. It is clear
that $c_u$ is a set of consecutive integers. Let $r(u)=\ds\max_{i\in
c_u}i$ and $l(u)=\ds\min_{i\in c_u}i$ denote the rightmost and the
leftmost cliques containing $u$ respectively. Note that two vertices
$u$ and $v$ are adjacent if and only if $c_u\cap c_v\ne\varnothing$.

Let $\eta:V\longrightarrow\integers$ be a labelling of vertices
obtained in the following manner: Choose a vertex $u_0$ such that
$r(u_0)\le r(v)$, $\forall v\ne u_0$. Assign label 0 to $u_0$ and all
vertices adjacent to $u_0$. Continue the same way considering only the
unlabelled vertices until all the vertices are labelled. More
formally:\vspace{5mm} 

\begin{algorithm}[H]\label{alg:labelVertex}
Let $V_0=V$, $I_\loc=\varnothing$, $i=0$\;
\While{$V_i\ne\varnothing$}{
$u_i\in V_i$ be such that $ r(u_i)\le r(v)$ $\forall v\in V_i$\;
$V'=\{u_i\}\cup\{v\in V_i| \textrm{$v$ is adjacent to $u_i$}\}$\;
$\eta(w)=i,\ \forall w\in V'$\;
$V_{i+1}=V_i\setminus V'$\;
$I_\loc\longleftarrow I_\loc\cup\{u_i\}$\;
$i\longleftarrow i+1$\;
}
\end{algorithm}

\begin{observation}\label{obs:lrrel}
For any vertex $v$, 
$\eta(v)\le i \Longleftrightarrow l(v)\le r(u_i)$.
\end{observation}
\begin{proof}
Since $v$ is adjacent to $u_{\eta(v)}$, we have $l(v)\le r(u_{\eta(v)})$.
It is clear that $r(u_{\eta(v)})\le r(u_i)$ since $\eta(v)\le
i$. Therefore, $l(v)\le r(u_i)$.

Suppose $\eta(v)> i$. From the algorithm, it implies that
$r(v)>r(u_i)$. Suppose $l(v)\le r(u_i)$, that is $l(v)\le r(u_i)\le
r(v)$. This implies that $v$ is adjacent to $u_i$. Then, by the algorithm
$\eta(v)\le i$, a contradiction. 
\qed
\end{proof}
\begin{observation}\label{obs:etauv}
For two vertices $v$ and $w$, if $\eta(v)=\eta(w)$, then $v$ and $w$
are adjacent.
\end{observation}
\begin{proof}
Let $\eta(v)=\eta(w)=i$. From Observation \ref{obs:lrrel} and from
the algorithm it follows that $l(v)\le r(u_i)\le r(v)$ and $l(w)\le
r(u_i)\le r(w)$. Therefore, $r(u_i)\in c_v\cap c_w$. Hence proved.
\end{proof}
In the algorithm let $l$ be the number of iterations, i.e.
$V_{l-1}\ne\varnothing$ and $V_l=\varnothing$.
\begin{observation}\label{obs:indset}
$I_\loc=\{u_0,u_1,\ldots,u_{l-1}\}$ is a maximum independent set. Hence,
$l=\alpha$.
\end{observation}
\begin{proof}
From the vertex numbering algorithm it is evident that $\loc$ is an
independent set. Suppose there exists an independent set of size greater
than $l$. By pigeon hole principle, at least two vertices in this set will
be assigned the same number and by Observation \ref{obs:etauv}, they will be
adjacent to each other, a contradiction.
\qed
\end{proof}
$I_\loc$ is crucial to our construction. From now on we refer to it as
 the \emph{primary independent set} with respect to the linear ordering
$\loc$.
\begin{observation}\label{obs:monotone}
$0=r(u_0)<r(u_1)<\cdots<r(u_{\alpha-1})=k-1$.
\end{observation}
\begin{proof}
From Observation \ref{obs:lrrel} we see that for
$i<\alpha-1$, $r(u_i)<l(u_{i+1})\le r(u_{i+1})$. Hence, 
$r(u_0)<r(u_1)<\cdots<r(u_{\alpha-1})$. Next we show that $r(u_0)=0$
and $r(u_{\alpha-1})=k-1$. 

Suppose, $r(u_0)\ne0$, then it is clear from the algorithm that for
all vertices $v$ with $l(v)=0$, $r(v)>0$. This implies that $C_0$
is a subset of $C_1$, which contradicts the maximality of the cliques.

It is easy to see that $r(u_{\alpha-1})\le k-1$. Suppose
$r(u_{\alpha-1})=t<k-1$. Consider any vertex $v\in C_{k-1}$. Clearly,
$r(v)=k-1>t$. Since $\eta(v)\le \alpha-1$, from Observation
\ref{obs:lrrel}, $l(v)\le t$. Therefore, $l(v)\le t\le r(v)$ which
implies $v\in C_t$. Hence, $C_{k-1}\subseteq C_t$, which contradicts
the maximality of the cliques.
\qed
\end{proof}
\subsection{Defining the Indifference Graphs}
Recall that $\loc:C_0, C_1,\ldots,C_{k-1}$ is a linear ordering of the
maximal cliques of $G$ and $I_\loc=\{u_0,\ldots,u_{\alpha-1}\}$ is
the primary independent set with respect to $\loc$. We can assume that
$\claw=2^p$, where $p$ is a positive integer. If not, we will work with
another interval graph $G'$ constructed in such a way that $\psi(G')=2^p$
 and $G$ is an induced subgraph of $G'$. To construct $G'$ from $G$
we consider a vertex $v\in C_{k-1}$. Let $m$ be the largest positive
integer such that there exists an induced $S(m)$ in $G$ with $v$ being
the central vertex of this $S(m)$. To obtain $G'$, we add $2^p-m$ new
vertices $v_{0},\ldots,v_{2^p-m-1}$ to $G$ such that they form an independent
set and are adjacent to only $v$. Then it is easy to verify that $G'$ would
correspond to the following linear ordering
of the maximal cliques: $\loc':C'_0,C'_1,\ldots,C'_{k+2^p-m-1}$,
where, $C'_i=C_i$ $0\le i\le k-1$ and $C'_{k+i}=\{v,v_i\}$ $0\le i\le
2^p-m-1$. Clearly, $\loc'$ satisfies Property \ref{prop:consec1} and
therefore $G'$ is an interval graph. Moreover,
we have an induced star $S(2^p)$ with $v$ as the central vertex. Clearly,
the remaining vertices of $G$ are unaffected by this construction. Hence,
$\psi(G')=2^p$. 

Now we define a function $f:\{0,\ldots,k-1\}\longrightarrow\reals$ as
follows:
\begin{enumerate}
\item $f(r(u_0))=f(0)=0$.
\item For $j\in\{ r(u_i)+1,\ldots,r(u_{i+1})\}$,
$f(j)=i+\frac{1}{2}+\frac{j- r(u_i)}{2(r(u_{i+1})- r(u_i))}$, for $0\le
i<\alpha-1$.
\end{enumerate}
\begin{remark}
From Observation \ref{obs:monotone} it is clear that $f$ is defined for
each $i\in\{0,1,\ldots,k-1\}$. Moreover,
$f$ is a strictly increasing function.
\end{remark}

Given positive integers $a$ and $i$, the \emph{$i$th bit function}
$b_i(\cdot)$ is defined as $b_i(a)=\floor{\frac{a}{2^i}}\mod 2$.
Now we define another labelling of vertices
$\gamma:V\longrightarrow\{0,1,\ldots,3\psi-1\}$ as follows:
\begin{equation}
\gamma(u)=\left\{
   \begin{array}{ll}
   \eta(u)\mod\psi + \psi, & \textrm{if $\floor{\frac{\eta(u)}{\psi}}$ is
   even},\vspace{2mm}\\
   \eta(u)\mod\psi + 2\psi, & \textrm{if $\floor{\frac{\eta(u)}{\psi}}$ is odd}.
   \end{array}\right.
\end{equation}
Recall that $p=\log_2\psi$. Note that $\gamma(u)$ is defined in such a
way that for $0\le i\le p-1$, $b_i(\gamma(u))=b_i(\eta(u))$, i.e. the
first $p$ bit positions of $\gamma(u)$ and $\eta(u)$ are identical. The
two extra bits in $p$th and $(p+1)$th positions depend on the parity
of $\floor{\frac{\eta(u)}{\psi}}$.

Now, we define $p+2=\log_2\psi+2$ indifference graphs $U_0, U_1,\ldots,U_{p+1}$ 
as follows. For each $U_i$ we define
$\Pi_i:V\longrightarrow\reals$ as per Definition \ref{defn:indiff}: For
$u\in V$,
\begin{equation}\label{eqn:assgn1}
\Pi_i(u)=\left\{
\begin{array}{ll}
f( r(u))-\psi+\frac{1}{2},& \textrm{if $b_i(\gamma(u))=0$},\\
f( l(u)),& \textrm{if $b_i(\gamma(u))=1$},
\end{array}\right. 
\end{equation}
where $0\le i\le p+1$. 
In the graph $U_i$, two vertices $u$ and $v$ are made adjacent
if and only if $|\Pi_i(v)-\Pi_i(u)|\le \psi-\frac{1}{2}$.

\subsection{Proof of $G=\ds\bigcap_{i=0}^{p+1}U_i$}
\begin{lemma}\label{lem:inclusion}
For any vertex $v$, $j\in c_v\Longrightarrow
f(j)\in\left[\Pi_i(v),\Pi_i(v)+\psi-\frac{1}{2}\right]$,
$0\le i\le p+1$.
\end{lemma}
\begin{proof}
Let $\eta(v)=m$. In order to handle some boundary cases, we define certain
notations. If $q<0$, then, let $r(u_q)=-1$. If $q>\alpha-1$,
then, let $r(u_q)=r(u_{\alpha-1})=k-1$.
\begin{clm}
$j\in c_v\Longrightarrow r(u_{m-1})+1\le j\le r(u_{m+\psi-1})$.
\end{clm}
\begin{proof}
If $m=0$, then it is clear that $l(v)=0=r(u_0)=r(u_{-1})+1$. Suppose $m>0$.
From Observation \ref{obs:lrrel} it immediately follows
that $l(v)\ge r(u_{m-1})+1$ and therefore $j>r(u_{m-1})$. 

Next, we show that $j\le r(u_{m+\psi-1})$.  Suppose
$m\ge\alpha-\psi$. Since $q=m+\psi-1\ge\alpha-1$, we have
$r(u_{m+\psi-1})=r(u_q)=r(u_{\alpha-1})=k-1$. But trivially, $j\le
 k-1$. Hence, we assume that $m<\alpha-\psi$. Suppose $v=u_m$, then
 this is trivially true from Observation \ref{obs:monotone}. Hence,
 we assume that $v\ne u_m$. Now, if there exists $j\in c_v$ such
 that $j>r(u_{m+\psi-1})$, then $t=r(u_{m+\psi-1})+1\in c_v$, since
by Observation \ref{obs:lrrel}, $l(v)\le r(u_{m+\psi-1})$ and $c_v$
is a set of consecutive integers. There exists a vertex $w\in C_t$ such
that $w\notin C_q$, for $q<t$, since otherwise $C_t$ will be a subset of
$C_{t-1}$. Clearly $w\ne v$. Now we claim that $\eta(w)=m+\psi$. Since
$l(w)=t>r(u_{m+\psi-1})$, by Observation \ref{obs:lrrel}, $\eta(w)\ge
m+\psi$. Also $l(u_{m+\psi})>r(u_{m+\psi-1})$ which implies
$r(u_{m+\psi})\ge l(u_{m+\psi})\ge t=l(w)$. By the algorithm, $r(w)\ge
 r(u_{m+\psi})$. Therefore, we have $l(w)\le r(u_{m+\psi})\le r(w)$
which implies that $w$ is adjacent to $u_{m+\psi}$, which in turn means
 $\eta(w)=m+\psi$. Since $v,w\in C_t$, they are adjacent. Clearly, the
vertex set $V'=\{u_m,u_{m+1},\ldots,u_{m+\psi-1},w\}$ forms an independent
 set since $l(w)=t>r(u_{m+\psi-1})$. Also, all the vertices of $V'$
 are adjacent to $v$ since, $l(v)\le r(u_m)\le r(u_{m+\psi-1})<l(w)\le
r(v)$. Therefore, $\{v\}\cup V'$ forms an induced star $S(\psi+1)$,
a contradiction. Hence, $j\le r(u_{m+\psi-1})$.
\bqed
\end{proof}
\begin{clm}\label{clm:psihalf}
$f(r(v))-f(l(v))<\psi-\frac{1}{2}$.
\end{clm}
\begin{proof}
From the above claim we have $r(u_{m-1})+1\le l(v)\le r(v)\le
r(u_{m+\psi-1})$. Now, by the definition of $f$ and noting
that $f$ is a strictly increasing function:
$\max\left(m-\frac{1}{2},0\right)<f(l(v))\le
f(r(v))\le\min\left(m+\psi-1,\alpha-1\right)$.
\bqed
\end{proof}
To complete the proof, we need to show that $\left[f(l(v)),f(r(v))\right]\subseteq
\left[\Pi_i(v),\Pi_i(v)+\psi-\frac{1}{2}\right]$.
If $b_i(\gamma(v))=0$,
\[
\left[\Pi_i(v),\Pi_i(v)+\psi-\frac{1}{2}\right]=\left[f(r(v))-
\psi+\frac{1}{2},f(r(v))\right],
\]
and if $b_i(\gamma(v))=1$,
\[
\left[\Pi_i(v),\Pi_i(v)+\psi-\frac{1}{2}\right]=\left[f(l(v)),
f(l(v))+\psi-\frac{1}{2}\right].
\]
In both cases it is sufficient to show that
 $f(l(v)>f(r(v))-\psi+\frac{1}{2}$, which immediately follows from
Claim \ref{clm:psihalf}.
\qed
\end{proof}

\begin{lemma}\label{lem:adj}
If $v,w\in V$ such that $v$ and $w$ are adjacent in $G$, then, $v$ and $w$
are adjacent in all the $p+2$ indifference graphs. 
\end{lemma}
\begin{proof}
Since $v$ and $w$ are adjacent, $c_v\cap c_w\ne\varnothing$. From Lemma
\ref{lem:inclusion} it follows that if $j\in c_v\cap c_w$, then, $f(j)\in
\left[\Pi_i(v), \Pi_i(v)+\psi-\frac{1}{2}\right] \cap
\left[\Pi_i(w),\Pi_i(w)+\psi-\frac{1}{2}\right]$ and hence,
$|\Pi_i(v)-\Pi_i(w)|\le \psi-\frac{1}{2}$ for $0\le i\le p+1$.
\qed
\end{proof}

\begin{lemma}\label{lem:nonadj}
If $v,w\in V$ such that $v$ and $w$ are not adjacent in $G$, then
there exists an indifference graph $U_i$, $i\in\{0,\ldots,p+1\}$, in which
$u$ and $w$ are not adjacent.
\end{lemma}
\begin{proof}
Without loss of generality we assume that $ r(v)< l(w)$. Since
$l(w)>r(v)\ge r(u_{\eta(v)})$, from Observation \ref{obs:lrrel} it
follows that $\eta(v)<\eta(w)$. 

Let $q_v=\floor{\frac{\eta(v)}{\psi}}$ and
$q_w=\floor{\frac{\eta(w)}{\psi}}$.  Now we consider
the following cases separately:
\begin{enumerate}
\item Suppose $q_w=q_v$: Then, $\gamma(v)\mod\psi<\gamma(w)\mod\psi$. This
in turn implies that there exists $i<\ceil{\log_2\psi}=p$ such that
$b_i(\gamma(v))=0$ and $b_i(\gamma(w))=1$. Then,
\begin{eqnarray*}
\Pi_i(w)-\Pi_i(v)=&&f( l(w))- f( r(v))+\psi-\frac{1}{2}>\psi-\frac{1}{2}.
\end{eqnarray*}
The last inequality follows from the fact that, by definition $f(\cdot)$ is a
strictly increasing function.
\item Suppose $q_w=q_v+1$: If $q_v$
is odd, then $b_{p}(\gamma(v))=0$ and $b_{p}(\gamma(w))=1$ and therefore,
as in Case 1, $\Pi_p(w)-\Pi_p(v)>\psi-\frac{1}{2}$. If
$q_v$ is even, then $b_{p+1}(\gamma(v))=0$ and
$b_{p+1}(\gamma(w))=1$ and similarly,
$\Pi_{p+1}(w)-\Pi_{p+1}(v)>\psi-\frac{1}{2}$. 
\item Suppose $q_w=q_v+2$: If $q_v$ is
even, then, $b_{p}(\gamma(v))=b_{p}(\gamma(w))=1$.
$\Pi_{p}(w)-\Pi_{p}(v)=f( l(w))-f( l(v))$. Note that $\eta(w)\ge q_w\psi$,
and therefore, from Observation \ref{obs:lrrel}, $l(w)\ge
r(u_{q_w\psi-1})+1$. Similarly, $\eta(v)\le q_v\psi+\psi-1$, and again from
Observation \ref{obs:lrrel}, $l(v)\le r(u_{q_v\psi+\psi-1})$. Therefore,  
\begin{eqnarray*}
f( l(w))-f( l(v)) &&\ge f(r(u_{q_w\psi-1})+1)-f(r(u_{q_v\psi+\psi-1}))\\ 
&&>\left(q_w\psi-1+\frac{1}{2}\right)-\left(q_v\psi+\psi-1\right)\\
&&=\left(q_v\psi+2\psi+\frac{1}{2}\right)-\left(q_v\psi+\psi\right)\\
&&=\psi+\frac{1}{2}>\psi-\frac{1}{2}.
\end{eqnarray*}
If $q_v$ is odd, then, $b_{p+1}(\gamma(v))=b_{p+1}(\gamma(w))=1$
and in a similar manner as above, we can show that
$\Pi_{p+1}(w)-\Pi_{p+1}(v)>\psi-\frac{1}{2}$.
\item Suppose $q_w>q_v+2$: If $b_{p}(\gamma(v))=b_{p}(\gamma(w))=1$, then,
we can show that $\Pi_{p}(w)-\Pi_{p}(v)>\psi-\frac{1}{2}$ in the same way
as Case 3. In a similar way, if
$b_{p+1}(\gamma(v))=b_{p+1}(\gamma(w))=1$, we can show that
$\Pi_{p+1}(w)-\Pi_{p+1}(v)>\psi-\frac{1}{2}$. 

Otherwise, from the definition of $\gamma(\cdot)$, it is easy to see
that either (1) $b_{p}(\gamma(v))=0$ and $b_{p}(\gamma(w))=1$ OR (2)
$b_{p+1}(\gamma(v))=0$ and $b_{p+1}(\gamma(w))=1$. As in Case 1 we
can show that $\Pi_{p}(w)-\Pi_{p}(v)>\psi-\frac{1}{2}$ for (1) and
$\Pi_{p+1}(w)-\Pi_{p+1}(v)>\psi-\frac{1}{2}$ for (2).
\end{enumerate}
Hence proved.
\qed
\end{proof}
Combining Lemmas \ref{lem:adj} and \ref{lem:nonadj}, we have
$G=\bigcap_{i=0}^{p+1}U_i$. Hence, we have proved Theorem
\ref{thm:mainThm}.

Note that when $\psi=\alpha$, the independence number of $G$, we have
 $q_v=0$ and therefore $b_p(\gamma(v))=1$ and $b_{p+1}(\gamma(v))=0$
for all vertices $v\in V$. From this, it is easy to see that $U_{p}$
and $U_{p+1}$ will correspond to complete graphs. Therefore, cubicity
of $G$ will be exactly $\clat$.

Next we observe that, given any interval graph $G$, we can construct a
graph $G'$ by adding a universal vertex to $G$. It is easy to that $G'$
 is an interval graph which contains $G$ as an induced subgraph. Also,
 $\psi(G')=\alpha(G')=\alpha$. By Lemma \ref{lem:inducedGraphCub},
it follows that $\cub(G)\le \cub(G')=\clat$. Considering this, Theorem
\ref{thm:mainThm} can be rewritten in the following way: \begin{theorem}
 Given an interval graph $G$, $\clct\le\cub(G)\le\min(\clct+2,\clat)$.
\end{theorem}

\begin{corollary}
Let $G$ be any graph. $\cub(G)\le\boxi(G)\clat$.
\end{corollary}
\begin{proof}
Let $b=\boxi(G)$. By Lemma \ref{lem:intcubbox}, there exist $b$ interval
graphs, say $G_i$, $0\le i<b$, such that $G=\bigcap_{i=0}^{b-1}G_i$. Since
each $G_i$ is a supergraph of $G$, $\alpha(G_i)\le \alpha$. Therefore,
$\cub(G_i)\le \clat$. Again by Lemma \ref{lem:intcubbox}, we have
$\cub(G)\le\sum_{i=0}^{b-1}\cub(G_i)\le\boxi(G)\clat$.
\qed
\end{proof}
We observe that this inequality is tight.   
In fact, given any two positive integers $k$ and $l$, there
exists a graph $G$ with $\boxi(G)=k$, $\alpha=l$ such that
$\cub(G)=k\ceil{\log_2l}$. One such example is the complete $k$-partite
graph with $|V|=kl$ (See Roberts \cite{recentProgressesInCombRoberts} for
proofs).


\end{document}